\documentclass{amsart}

\usepackage{lineno,hyperref}

\usepackage{latexsym,amsmath,amsthm}
\usepackage{amssymb}

\def\ds{\displaystyle}

\newcommand{\R}{{\mathbb R}}
\newcommand{\m}{\mathfrak{m}}
\def\ba{{\mathbf a}}
\def\bs{\boldsymbol}

\newtheorem{thm}{Theorem}[section]  
\newtheorem{lem}[thm]{Lemma}	       

\numberwithin{equation}{section}

\begin{document}

\title[Non-homogeneous conormal derivative problems]{Non-homogeneous conormal derivative problems for discontinuous quasilinear elliptic equations with Morrey data}

\author[D.K. Palagachev]{Dian K. Palagachev}
\address{Department of Mechanics, Mathematics and Management,
Politechnic University of Bary, Italy}
\email{ dian.palagachev@poliba.it}

\author[L.G. Softova]{Lubomira G. Softova}
\address{Department of Mathematics,
University of Salerno, Italy}
\email{ lsoftova@unisa.it}



\date{\today}

\begin{abstract}
A non-homogeneous conormal derivative problem is considered for quasilinear di\-ver\-gen\-ce form elliptic equations modeled on the $m$-Laplacian operator. The nonlinear
terms are given by Carath\'eodory functions and satisfy controlled growth structure conditions
with respect to the solution and its gradient,
while their $x$-behaviour is controlled in terms of suitable Morrey spaces.

Global essential boundedness is proved
for the weak solutions, generalizing thus the classical $L^p$-result of Ladyzhenskaya and Ural'tseva to the framework of the Morrey  scales.
\end{abstract}

\maketitle

\section{Introduction}

The main goal of the present paper is to prove essential boundedness of the $W^{1,m}(\Omega)$-weak solutions, $m\in (1,n],$ of the non-homogeneous conormal derivative problem for second-order divergence form elliptic equations
\begin{equation}\label{1.1}
\begin{cases}
\ds\text{div}\, \ba(x,u,Du)= b(x,u,Du) 
& x\in \Omega, \\
\ds\ba(x,u,Du)\cdot \bs{\nu}(x) =\psi(x,u) \quad  & x\in  \partial\Omega.
\end{cases}
\end{equation}
Here $\Omega \in \R^n,$ $n \geq 2,$ is a bounded domain with Lipschitz continuous boundary, and $\bs{\nu}(x)=\big(\nu_1(x),\ldots,\nu_n(x)\big)$ denotes the unit outer normal to $\partial\Omega,$ and the nonlinear terms 
$\ba(x,z,\xi)=\big(a_1(x,z,\xi),\ldots,a_n(x,z,\xi)\big),$ $b(x,z,\xi)$  and $\psi(x,z)$ are  of Carath\'eodory type, that is, these are measurable in $x\in \Omega$ (or in $x\in\partial\Omega$ in the case of $\psi$) for all    $z\in \R,$ $\xi \in \R^n$ and continuous with respect to $(z,\xi)$ for almost all $x\in \Omega/\partial\Omega.$ 

The divergence form nonlinear operator is modelled after the $m$-Laplacian operator with an arbitrary exponent $m\in(1,n],$ and thus 
the following \textit{coercivity condition} of order $m$ is assumed in the sequel:
there exist positive constants $\gamma$ and $\Lambda$ such that 
\begin{equation}\label{1.3} 
    \ba(x,z,\xi)\cdot \xi\geq \gamma |\xi|^m
    -\Lambda|z|^{{m^*}}-\Lambda\varphi_1(x)^\frac{m}{m-1}
\end{equation}
for almost all (a.a.) $x \in \Omega$ and  all $(z,\xi) \in \R\times \R^n,$ where $\varphi_1\in L^{\frac{m}{m-1}}(\Omega)$ and 
 $m^*$ denotes the Sobolev conjugate of $m,$ that is, $m^*=\frac{nm}{n-m}$ if $m<n$ and $m^*$ is any exponent greater than $n$ when $m=n.$

Besides  the \textit{coercivity} condition \eqref{1.3}, we assume that the nonlinear terms in \eqref{1.1} support \textit{controlled growths} with respect to $u$ and $|Du|,$ that is, 
\begin{equation}\label{1.4}
|\ba(x,z,\xi)|\leq  \Lambda \Big(\varphi_1(x)  +|z|^{\frac{m^*(m-1)}{m}}+|\xi|^{m-1}\Big) 
\end{equation}
with $\varphi_1$ as in \eqref{1.3},
\begin{equation}\label{1.6}
|b(x,z,\xi)|\leq \Lambda \Big(\varphi_2(x)+|z|^{m^*-1} +|\xi|^{\frac{m(m^*-1)}{m^*}}\Big)
\end{equation}
with $\varphi_2\in L^{\frac{mn}{mn+m-n}}(\Omega)$ 
and
\begin{equation}\label{1.8}
|\psi(x,z)|\leq \psi_1(x)+\psi_2(x)|z|^\beta, 
\end{equation}
where $\beta\in\left[0,\frac{n(m-1)}{n-m}\right),$  $\psi_1\in L^{\frac{m(n-1)}{n(m-1)}}(\partial\Omega),$
$\psi_2\in L^{\frac{m(n-1)}{n(m-1)-\beta(n-m)}}(\partial\Omega)$ and $\frac{n(m-1)}{n-m}$ is to be intended as $+\infty$ in case $m=n.$

\medskip

Recall that, fixed a real number $m\in(1,n],$ a function $u \in W^{1,m}(\Omega)$ is  a \textit{weak solution} of the problem \eqref{1.1} when the standard integral identity
\begin{align}\label{1.2}
\int_\Omega& \ba\big(x,u(x),Du(x)\big)\cdot Dv(x)\,dx+\int_\Omega b\big(x,u(x),Du(x)\big)v(x) \, dx\\
\nonumber
     &\qquad\qquad \qquad=\int_{\partial\Omega} \psi\big(x,u\big)v(x)\, d\Gamma_x
\end{align}
is satsfied for each  $v \in W^{1,m}(\Omega).$ Actually, the integrals here involved
converge because of the controlled growths conditions \eqref{1.4}, \eqref{1.6} and \eqref{1.8} and the required integrability of the functions governing the $x$-behaviours therein, and the controlled growths are the minimal assumptions guaranteeing sense of  the concept of weak solution to \eqref{1.1}.

\medskip

In the case of \textit{sub-controlled} terms   $\ba(x,z,\xi)$ and $b(x,z,\xi)$ when these grow as $|z|^{m-1}$ and $|\xi|^{m-1},$ 
boundedness of the classical solutions  $(u\in C^1(\overline{\Omega})\cap C^2(\Omega))$ to quasilinear conormal problems has been studied by Lieberman in \cite{Lieb1,Lieb2} (see also \cite{Lieb3}). In the same situation, Winkert derived in \cite{Win} an $L^\infty$-\textit{a~priori} bound for the weak solutions of \eqref{1.1}, and similar result has been proved in \cite{WinZach} assuming variable growth exponents in the nonlinear terms of \eqref{1.1}. For what concerns the case of \textit{controlled growths} of the nonlinearities,  there is a classical result due to Ladyzhenskaya and Ural'tseva \cite[Chapter~X, \S~2]{LU} which asserts boundedness and H\"older continuity of the weak solutions when the  $x$-behaviour of $\ba(x,z,\xi)$ and $b(x,z,\xi)$ is controlled in terms of suitable Lebesgue spaces (see also \cite{Kim} in the particular case when $m=2$). Recently, boundedness and H\"older continuity have been obtained in \cite{H-K-Win-Z} for the non-homogeneous problem \eqref{1.1} in a very general \textit{controlled growths} r\'egime when the nonlinearities are subject to variable exponents growths. Anyway, these papers do not allow singularities with respect to $x$ of the  nonlinear terms in \eqref{1.1}, restricting their $x$-behaviour to $L^\infty.$
We have to note also the deep papers of Arkhipova \cite{Arkh1,Arkh2}, where reverse H\"older inequalities have been derived for solutions of quasilinear conormal derivative problems which lead to \textit{improving of gradient integrability} and the subsequent  regularity results. 

Having in mind that the
\textit{controlled growths} are \textit{optimal} (cf. the counterexamples in \cite[Chapter~I, \S~2]{LU}) in order to have  boundedness of the weak solutions, our goal here is to weaken the hypotheses on the functions $\varphi_1,\varphi_2,\psi_1$ and $\psi_2$ that govern the $x$-behaviour of $\ba(x,z,\xi)$ and $b(x,z,\xi).$ Namely, we will take these in \textit{Morrey spaces}  with suitable exponents (see Section~\ref{sec2} for the definition of the Morrey spaces). Thus,
 regarding the behaviour of the nonlinear terms with respect to $x,$ we suppose 
\begin{equation}\label{1.5}
\begin{cases}
\varphi_1\in L^{p_1,\lambda_1}(\Omega),\quad & \lambda_1\in [0,n),\quad
p_1>\max\left\{\dfrac{m}{m-1},\dfrac{n-\lambda_1}{m-1}\right\},\\[10pt]
\varphi_2\in L^{p_2,\lambda_2}(\Omega),\quad & \lambda_2\in [0,n),\quad
p_2>\max\left\{\dfrac{mn}{mn-n+m},\dfrac{n-\lambda_2}{m}\right\}
\end{cases}
\end{equation}
and 
\begin{equation}\label{Q}
\begin{cases}
\psi_1\in L^{q_1,\mu_1}(\partial\Omega),\quad & \mu_1\in [0,n-1),\\
\quad &
q_1>\max\left\{\dfrac{m(n-1)}{n(m-1)},\dfrac{n-1-\mu_1}{m-1}\right\},\\[10pt]
\psi_2\in L^{q_2,\mu_2}(\partial\Omega),\quad & \mu_2\in [0,n-1),\\
\quad & 
q_2>\max\left\{\dfrac{m(n-1)}{n(m-1)-\beta(n-m)},\dfrac{n-1-\mu_2}{m-1}\right\}
\end{cases}
\end{equation}
with $\beta\in\left[0,\frac{n(m-1)}{n-m}\right)$ as before.

Recalling that $\varphi_1\in L^{\frac{m}{m-1}}(\Omega),$  $\varphi_2\in L^{\frac{mn}{mn+m-n}}(\Omega),$ $\psi_1\in L^{\frac{m(n-1)}{n(m-1)}}(\partial\Omega)$ and 
$\psi_2\in L^{\frac{m(n-1)}{n(m-1)-\beta(n-m)}}(\partial\Omega)$
are necessary conditions, guaranteeing convergence of the integrals  in \eqref{1.2}, we need here the  stronger  hypotheses $p_1>\frac{m}{m-1},$  $p_2>\frac{mn}{mn-n+m},$ 
$q_1>\frac{m(n-1)}{n(m-1)}$ and
$q_2>\frac{m(n-1)}{n(m-1)-\beta(n-m)}$
in order to gain \textit{better integrability of the gradient}  through the results of Arkhipova \cite{Arkh1,Arkh2}.

\medskip

The main result of the paper (Theorem~\ref{thm3.1}) claims \textit{essential boundedness} of any $W^{1,m}(\Omega)$-weak solution of the conormal derivative problem \eqref{1.1} under the above mentioned assumptions. Moreover, the $L^\infty(\Omega)$-norm of the solution turns out to be bounded in terms of the data involved in the hypotheses imposed, and also of $\|u\|_{W^{1,m}(\Omega)}.$ Let us stress the reader attention at this point that the dependence of boundedness on the $W^{1,m}(\Omega)$-norm of the solution itself is \textit{not restrictive} because our result is \textit{qualitative} and no \textit{quantitative} as the \textit{$L^\infty$-a~priori} estimates provided in \cite{Win,WinZach} for instance. The price for this is paid by the possibility to deal with \textit{strictly controlled growths} of the nonlinearities -- a situation when additional restrictions must be imposed if one wants to obtain a real \textit{$L^\infty$-a~priori} estimate. Let us note that boundedness of the weak solutions to \eqref{1.1} has already been proved in our previous paper \cite{AFPS} under the assumption
\begin{equation}\label{hbc}
\psi(x,z)z\leq0\quad \text{for a.a.}\ x\in\partial\Omega\ \text{and all}\ z\in\R.
\end{equation}
This \textit{sign-type} requirement on the right-hand side of the boundary condition is so restrictive that it actually makes disappear the surface integral in \eqref{1.2}, exactly as happens if \eqref{1.1} would be a \textit{homogeneous} problem, that is, when $\psi\equiv0.$ On the other hand, the evaluation of the surface integral term in \eqref{1.2} offers new and interesting chalenges related to Sobolev inequalities of trace type and other measure density results which, at the end, allow to study a really \textit{non-homogeneous} conormal derivative problems. 

Indeed, taking $\lambda_1=\lambda_2=0$ in \eqref{1.5} and assuming the sing-condition \eqref{hbc}, our Theorem~\ref{thm3.1} recovers the boundedness result of Ladyzhenskaya and Ural'tseva (\cite[Chapter~X, \S~2]{LU}), and  the restrictions $p_1>\frac{n}{m-1},$ $p_2>\frac{n}{m}$ and $q_1,q_2>\frac{n-1}{m-1}$ 
are sharp when working in the framework of the Lebesgue spaces as known by the counterexamples in \cite[Chapter~I, \S~2]{LU}. Our boundedness result shows that, taking $\varphi_1,$  $\varphi_2,$ $\psi_1$ and $\psi_2$ in Morrey spaces, the values of $p_1,$ $p_2,$ $q_1$ and $q_2$ could be even decreased  at the expense of increase $\lambda_1,$ $\lambda_2,$ $\mu_1$ and $\mu_2,$ respectively, still maintaining the restrictions $p_1>\frac{n-\lambda_1}{m-1},$ $p_2>\frac{n-\lambda_2}{m},$  $q_1>\frac{n-1-\mu_1}{m-1}$ and $q_2>\frac{n-1-\mu_2}{m-1}.$

It is interesting to note that further regularity properties of the weak solutions to conormal derivative problems could be obtained by combining our global boundedness result with the techniques developed by Nazarov and Ural'tseva in \cite{NU} where local properties, such as strong maximum principle, Harnack inequality and H\"older continuity, have been obtained for weak solutions to \textit{linear} divergence form equations with Morrey lower-order coefficients.

\medskip

The machinery used in the proof of Theorem~\ref{thm3.1} is that of \cite{BPS} (where quasilinear Dirichlet problems were considered) and \cite{AFPS}, and it relies on the
De Giorgi approach to the boundedness as adapted by Ladyzhenskaya and Ural'tseva
(cf. \cite[Chapter~IV]{LU}) to the quasilinear situation. Precisely, the controlled growth
assumptions allow to obtain nice decay estimates for the total mass of the weak solution taken over its level sets. However, unlike the $L^p$-approach of Ladyzhenskaya and
Ural'tseva, the mass we have to do with is taken with respect to a suitable positive Radon
measure $\m$ that depends not only on the Lebesgue measure, but also on $\varphi_1^{\frac{m}{m-1}},$ $\varphi_2$
  and a suitable power of the solution $u$ itself. Thanks to the Morrey integrability of the data (the hypotheses \eqref{1.5}),
precise trace-type Sobolev inequalities allow to estimate the $\m$-mass of $u$ in terms of the $m$-energy of $u.$ At this point  the controlled growths assumptions combine with
the better-gradient-integrability results of Arkhipova (\cite{Arkh1,Arkh2}) in order to estimate the $m$-energy of $u$ in
terms of a small multiplier of the same quantity plus a suitable power of the $\m$-measure of the solution level
set. The global boundedness of the weak solution then follows by a
classical result known as \textit{Hartman--Stampacchia maximum principle} that ensures \textit{finite time extinction property} of  non-increasing functions with suitable decay at infinity.

\medskip

In what follows, we will use the omnibus phrase 
``\textit{known quantities}'', meaning that a given constant depends on the data in all the above presented hypotheses,  which include  $n,$ $m,$ 
 $\gamma,$ $\Lambda,$ $p_1,$ $p_2,$ $\lambda_1,$ $\lambda_2,$ $q_1,$ $q_2,$ $\mu_1,$ $\mu_2,$ $\beta,$  $\|\varphi_1\|_{L^{p_1,\lambda_1}(\Omega)},$ $\|\varphi_2\|_{L^{p_2,\lambda_2}(\Omega)},$ $\|\psi_1\|_{L^{q_1,\mu_1}(\partial\Omega)},$ $\|\psi_2\|_{L^{q_2,\mu_2}(\partial\Omega)},$
 $\mathrm{diam}\,\Omega$ and the Lipschitz regularity of $\partial\Omega.$ The letter $C$ will denote a generic constant, depending on known quantities, which may vary within the same formula. As already said, we are going to obtain a \textit{qualitative} boundedness result and not an \textit{a~priori} estimate. In this sense, the dependence of certain constants also on the $W^{1,m}(\Omega)$-norm of the weak solution is not restrictive, and we will consider in the sequel also 
$\|u\|_{W^{1,m}(\Omega)}$ as a \textit{known quantity}. When necessary, the dependence on $\|u\|_{W^{1,m}(\Omega)}$ will be explicitly mentioned.

\section{Functional analysis tools}\label{sec2}
 
For readers convenience we collect here some basic   functional analysis results to be used in the sequel. In what follows, $\Omega\subset\R^n$ will be a bounded domain with Lipschitz continuous boundary.

\medskip 

The Morrey space $L^{p,\lambda}(\Omega),$
 $p\in[1,\infty)$ and $\lambda\in[0,n],$  consists of all functions $\varphi\in L^p(\Omega)$ such that 
$$
\|\varphi\|_{L^{p,\lambda}(\Omega)}:=\sup_{x_0\in\Omega,\ \rho>0}\left(\rho^{-\lambda}\int_{B_\rho(x_0)\cap\Omega}|\varphi(x)|^p\, dx\right)^{1/p} < \infty,
$$
where $B_\rho(x_0)$ stands for a ball of radius $\rho$ centered at $x_0.$ 
Under the norm given by the last quantity, 
$L^{p,\lambda}(\Omega)$  becomes a Banach space, and the limit cases $\lambda=0$ and $\lambda=n$ correspond to $L^p(\Omega)$ and $L^\infty(\Omega),$ respectively.

\begin{lem}\label{lem2.1}
{\em (Embeddings between Morrey spaces, cf. \cite{Pic})}
For arbitrary $p,\bar{p}\in [1,\infty)$ and $\lambda,\bar{\lambda} \in [0,n),$ the inclusion
$$
L^{p,\lambda}(\Omega)\subseteq L^{\bar{p},\bar{\lambda}}(\Omega)
$$
holds \em{if and only if}
$
p\geq \bar{p}\geq 1$ and $\dfrac{p}{n-\lambda}\geq\dfrac{\bar{p}}{n-\bar{\lambda}}.$
\end{lem}

The Morrey spaces on the boundary $\partial\Omega$ are defined in a similar manner employing, however, the natural localization and flattening (cf. \cite{Cam}). Thus, bearing in mind the Lipschitz continuity of $\partial\Omega,$ there is a positive $R$ such that for each point $x_0\in\partial\Omega$ there is an one-to-one Lipschitz mapping $\Psi$ of  $B_{R}(x_0)$ \textit{onto} $D\subset\R^n$ with the property that $\Psi(B_{R}(x_0)\cap \Omega)\subset \R^n_+=\{y=(y',y_n)\in\R^n\colon\ y_n>0\}$
and $\Psi(B_{R}(x_0)\cap \partial\Omega)=D'\subset \partial\R^n_+.$ Given a function $\psi\in L^q(\partial\Omega),$  $q\geq1$ and $\mu\in[0,n-1],$ we set $\psi':=\psi\circ\Psi^{-1}$ and define
$$
\|\psi'\|_{L^{q,\mu}(D')}:=\sup_{x'_0\in D',\ \rho>0}\left(\rho^{-\mu}\int_{B'_\rho(x'_0)\cap D'}|\psi(\Psi^{-1}(y'))|^q\, dy'\right)^{1/q} 
$$
where $B'_\rho(x'_0)$ is the $(n-1)$-dimensional  ball of radius $\rho$ centered at $x'_0.$ Since $\partial\Omega$ is a compact, it can be covered by a finite number of balls $B_R(x_k),$ $x_k\in\partial\Omega,$ $k=1,\ldots,N,$ and to any $\psi\in L^q(\partial\Omega)$ there correspond $k$ domains $D_k$ and the respective diffeomorphisms $\Psi_k$ that define $k$ functions $\psi'_k\colon D'_k\to\R,$ given by $\psi\circ\Psi_k^{-1}.$

This way, the Morrey space $L^{q,\mu}(\partial\Omega)$ is the collection of all functions $\psi\in L^q(\partial\Omega)$ for which the quantities $\|\psi'_k\|_{L^{q,\mu}(D'_k)}$ are finite, and we set
$$
\|\psi\|_{L^{q,\mu}(\partial\Omega)}:= \sum_{k=1}^n \|\psi'_k\|_{L^{q,\mu}(D'_k)}.
$$

\bigskip

In what follows, we will often use the following \textit{multiplicative} version of the \textit{Gagliardo--Nirenberg interpolation inequality,} where $\langle\langle v\rangle\rangle_\sigma:=\left(\int_\Omega |v(x)|^\sigma\,dx\right)^{1/\sigma}$ for $\sigma>0.$
\begin{lem}\label{G-N}
{\em (cf. \cite[Theorem 1.4.8/1]{Maz})}
Let $\Omega$ be a bounded domain with the cone property and $v\in W^{1,m}(\Omega),$ $m\in[1,n).$

Then there is a constant $C$ such that
$$
\langle\langle v\rangle\rangle_r \leq C \Big( \langle\langle Dv\rangle\rangle_m + \langle\langle v\rangle\rangle_\sigma\Big)^\delta \langle\langle v\rangle\rangle_\sigma^{1-\delta}
$$
for all $\delta\in[0,1]$ and all $\sigma\in\left(0,\dfrac{nm}{n-m}\right],$ where $\dfrac{1}{r}=\delta\dfrac{n-m}{nm}+\dfrac{1-\delta}{\sigma}.$
\end{lem}

The next \textit{trace}-type theorem, rewritten in our situation, will play a crucial r\^ole in the estimate of the surface integral in \eqref{1.2}:
\begin{lem}\label{trace}
{\em (cf. \cite[Corollary 1.4.7/2]{Maz})}
Let $\Omega\subset\R^n$ be a bounded and Lipschitz domain  and let $\m$ be a measure in $\overline{\Omega}$ with the property
\begin{equation}\label{measure}
\sup_{x\in\overline{\Omega},\ \rho\in(0,1)} \rho^{-s}\m(B_\rho(x))\leq K 
\end{equation}
for some $s\in(n-m,n]$ with a constant $K.$ 

Then there is an absolute constant $C=C(n,s,m,K,\partial\Omega)$ such that
$$
\|v\|_{L^r(\m,\overline{\Omega})}:=
\left(\int_\Omega |v(x)|^r\,d\m\right)^{1/r}
\leq C \|v\|_{W^{1,m}(\Omega)}^\tau \|v\|_{L^{m}(\Omega)}^{1-\tau}
$$
for each $v\in W^{1,m}(\Omega),$ $m\in[1,n),$ where $r\in\left[m,\dfrac{sm}{n-m}\right)$ and $\tau=\dfrac{n}{m}-\dfrac{s}{r}<1.$
\end{lem}

The first result of this kind, in a \textit{non-multiplicative} version, has been proved by D.R.~Adams in \cite{Ad,Ad'} with regard to Riesz potentials and functions in the Sobolev space $W^{1,m}_0(\Omega).$ Precisely, under the hypothesis \eqref{measure}, one has
$$
\|v\|_{L^\frac{sm}{n-m}(\m,\overline{\Omega})}\leq C \|Dv\|_{L^m(\Omega)}\qquad \forall v\in W^{1,m}_0(\Omega).
$$
(Note that this is exactly what Lemma~\ref{trace} asserts with $\tau=1$ formally!) Willing to extend that inequality to functions $v\in W^{1,m}(\Omega)$ with non necessarily zero boundary trace,  it suffices to extend $v$ to $V\in W^{1,m}(\R^n)$ in a way that $\|V\|_{W^{1,m}(\R^n)}\leq C(\partial\Omega)\|v\|_{W^{1,m}(\Omega)}$ and then multiply by a suitable cut-off function $\zeta$  such that $\zeta\equiv 1$ over $\Omega.$ Assuming that the measure $\m$ is extended as zero outside $\Omega,$ application of the original Adams trace inequality to $\zeta V$ then yields
$$
\|v\|_{L^\frac{sm}{n-m}(\m,\overline{\Omega})}\leq C \left(\int_\Omega\big(|Dv(x)|^m+|v(x)|^m\big)\, dx\right)^{1/m}\qquad \forall v\in W^{1,m}(\Omega),
$$
and interpolation of $\int_\Omega |v(x)|^m\,dx$ by using Lemma~\ref{G-N} and Young's inequality implies
\begin{lem}\label{lem2.2}
{\em (Adams trace inequality)}
Let $\m$ be a  measure in $\overline{\Omega}$  satisfying \eqref{measure}.

Then there exists a constant $C=C(n,s,m,\sigma,K,\partial\Omega)$ such that
$$
\|v\|_{L^\frac{sm}{n-m}(\m,\overline{\Omega})}\leq C\left(\left(\int_\Omega|Dv(x)|^m\, dx\right)^{1/m}+\left(\int_\Omega|v(x)|^\sigma\, dx\right)^{1/\sigma}\right)
$$
for all $v\in W^{1,m}(\Omega)$ and all $\sigma\in (0,m].$

In particular, if $d\m=c(x)\,dx$ with $c\in L^{1,n-m+\varepsilon_0}(\Omega)$ and small $\varepsilon_0>0,$ then there is a constant $C=C\left(n,\varepsilon_0,m,\sigma,\|c\|_{L^{1,n-m+\varepsilon_0}(\Omega)},\partial\Omega\right)$ such that
\begin{align*}
& \left(\int_\Omega|v(x)|^\frac{(n-m+\varepsilon_0)m}{n-m} c(x)\, dx\right)^\frac{n-m}{(n-m+\varepsilon_0)m}\\
&\qquad\qquad \leq C
\left(\left(\int_\Omega|Dv(x)|^m\, dx\right)^{1/m}+\left(\int_\Omega|v(x)|^\sigma\, dx\right)^{1/\sigma}\right)
\end{align*}
for all  $v\in W^{1,m}(\Omega)$ and all $\sigma\in (0,m].$
\end{lem}

In the particular case when $\m$ is the Lebesgue measure supported in $\Omega,$ we have $s=n$ in \eqref{measure}  and Lemma~\ref{lem2.2} implies the next version of the Sobolev embedding theorem, valid for functions with non necessary zero boundary trace.
\begin{lem}\label{lem2.3}
{\em (Sobolev inequality, see also \cite[Chapter~II, \S~2]{LU})} 
Let $\Omega\subset\R^n$ be a bounded domain with Lipschitz continuous boundary and $m<n.$
 
Then there exists a constant $C$ depending on $n,m,$ $\text{\em diam}\,\Omega$ and the Lipschitz regularity of $\partial\Omega,$ such that
$$
\|v\|_{L^{m^*}(\Omega)}\leq C\left(\left(\int_\Omega|Dv(x)|^m\, dx\right)^{1/m}+\left(\int_\Omega|v(x)|^\sigma\, dx\right)^{1/\sigma}\right)
$$
for all $v\in W^{1,m}(\Omega)$ and all $\sigma\in (0,m].$
\end{lem}

We conclude the list of auxiliary results with the next one, known as \textit{Hartman--Stampacchia maximum principle}, which asserts \textit{finite time extinction property} of the non-increasing functions with suitable decay at infinity:
\begin{lem}\label{lem2.5}
{\em (cf. \cite{HS}, \cite[Chapter~II, Lemma~5.1]{LU})}
Let $\tau\colon\R\to[0,\infty)$ be a  non-increasing function and suppose there exist constants $C>0,$ $k_0\geq0,$ $\delta>0$ and 
$\alpha \in [0,1+\delta]$ such that
$$
\int_k^\infty \tau(t)\, dt \leq Ck^\alpha\big(\tau(k)\big)^{1+\delta} \quad
\forall k\geq k_0.
$$
Then there is a number $k_{\max},$ depending on $C,$ $k_0,$ $\delta,$ $\alpha$ and $\int_{k_0}^\infty \tau(t)\, dt,$ such that
$$
\tau(k)=0\quad \forall k\geq k_{\max}.
$$
\end{lem}

\section{The Main Result}

The main result of the present paper asserts global essential boundedness of any weak solution to the non-homogeneous conormal problem \eqref{1.1}:
\begin{thm}\label{thm3.1}
Let $\Omega \in \R^n$ be a bounded domain with Lipschitz continuous boundary, $m\in(1,n],$ and suppose the hypotheses \eqref{1.3}--\eqref{1.8} together with \eqref{1.5} and \eqref{Q} are satisfied.

Then any $W^{1,m}(\Omega)$-weak solution of the non-homogeneous conormal problem \eqref{1.1} is globally essentially bounded. Precisely, 
\begin{equation}\label{3.1}
    \| u\|_{L^{\infty}(\Omega)} \leq M
\end{equation}
with a constant depending on known quantities, on $\|u\|_{L^m(\Omega)}$ and on $\| Du\|_{L^m(\Omega)}.$ 
\end{thm}

A crucial step in the \textsc{PROOF} of Theorem~\ref{thm3.1}, which allows the treatment of the borderline values of the exponents in the $u,|Du|$-controlled growths conditions, is ensured by the \textit{higher gradient integrability} property of the weak solutions to \eqref{1.1}. 
This \textit{Gehring--Giaquinta} type result requires only coercivity and controlled growths of the nonlinear ingredients in \eqref{1.1}, together with some more integrability of the functions, governing their $x$-behaviour. The reader is refered to \cite[Theorem~2]{Arkh1} and \cite[Theorem~4]{Arkh2} with the subsequent remarks for the corresponding proof (see also \cite[Section~8]{Kim} for the particular case  $m=2$).
\begin{lem}\label{lem2.4}
Assume \eqref{1.3}--\eqref{1.8} together with $\varphi_1\in L^{p_1}(\Omega),$ $p_1>\frac{m}{m-1},$ $\varphi_2\in L^{p_2}(\Omega),$ $p_2>\frac{mn}{mn-n+m},$ $\psi_1\in L^{q_1}(\partial\Omega),$ $q_1>\frac{m(n-1)}{n(m-1)}$ and
$\psi_2\in L^{q_2}(\partial\Omega),$ $q_2>\frac{m(n-1)}{n(m-1)-\beta(n-m)}$ with 
$\beta\in\left[0,\frac{n(m-1)}{n-m}\right).$

Then there is an exponent $m_0>m$ such that
any weak solution $u\in W^{1,m}(\Omega),$ $m\in(1,n],$ of the non-homogeneous conormal derivative problem \eqref{1.1} belongs to
 $ W^{1,m_0}(\Omega),$ and 
$$
\|Du\|_{L^{m_0}(\Omega)}\leq N,
$$
with a constant $N$ depending on known quantities and on  $\|Du\|_{L^{m}(\Omega)}.$
\end{lem}

Starting with the \textsc{PROOF} of Theorem~\ref{thm3.1}, let us assume 
$m<n$ firstly and define the measure
$$
d\m:=\left(\chi_\Omega(x)+\varphi_1(x)^\frac{m}{m-1}+\varphi_2(x)+|u(x)|^\frac{m^2}{n-m}\right)\,dx, 
$$
where $\chi_\Omega(x)$ is the indicator function of $\Omega$  and, with no loss of generality, $\varphi_1$ and $\varphi_2$ are supposed to be extended as zero outside $\Omega.$

We want  to check now that $\m$ is a measure for which the Adams trace inequality (Lemma~\ref{lem2.2}) is applicable. For, 
it follows from the Lipschitz continuity of $\partial\Omega$ and \eqref{1.5} that for any
ball  $B_\rho$ of radius $\rho$ one has
\begin{align*}
\int_{B_\rho}\chi_\Omega(x)\, dx &= \int_{B_\rho\cap\Omega}\, dx \leq C \rho^n = C\rho^{n-m+m},\\
\int_{B_\rho}\varphi_1(x)^\frac{m}{m-1}\, dx &\leq \|\varphi_1\|_{L^{p_1,\lambda_1}(\Omega)}^\frac{m}{m-1}\rho^{n-\frac{m(n-\lambda_1)}{p_1(m-1)}}
=\|\varphi_1\|^\frac{m}{m-1}_{L^{p_1,\lambda_1}(\Omega)}\rho^{n-m+\left(m-\frac{m(n-\lambda_1)}{p_1(m-1)}\right)}, \\
\int_{B_\rho}\varphi_2(x)\, dx &\leq \|\varphi_2\|_{L^{p_2,\lambda_2}(\Omega)}\rho^{n-\frac{n-\lambda_2}{p_2}} = \|\varphi_2\|_{L^{p_2,\lambda_2}(\Omega)} \rho^{n-m+\left(m-\frac{n-\lambda_2}{p_2}\right)}
\end{align*}
with $m-\frac{m(n-\lambda_1)}{p_1(m-1)}>0$ and $m-\frac{n-\lambda_2}{p_2}>0$ because of the assumptions $p_1>\frac{n-\lambda_1}{m-1}$ and $p_2>\frac{n-\lambda_2}{m}.$ Further on, bearing in mind that $u\in L^{m_0^*}(\Omega)$ as follows from Lemmas~\ref{lem2.4} and \ref{lem2.3},  the H\"older inequality implies 
\begin{equation}\label{3.2}
\int_{B_\rho}|u(x)|^\frac{m^2}{n-m}\, dx \leq \|u\|^\frac{m^2}{n-m}_{L^{m^*_0}(\Omega)} \rho^{n-m+\left(m-\frac{nm^2}{m^*_0(n-m)}\right)}
\end{equation}
with $m-\frac{nm^2}{m^*_0(n-m)}>0$ since $m^*_0>{m^*}=\frac{nm}{n-m}.$

Setting
$$
\varepsilon_0:=\min\left\{m,m-\frac{m(n-\lambda_1)}{p_1(m-1)},m-\frac{n-\lambda_2}{p_2},m-\frac{nm^2}{m^*_0(n-m)}\right\}>0,
$$
we get
$$
\m(B_\rho)\leq K \rho^{n-m+\varepsilon_0}
$$
with a constant $K$ depending on known quantities, and on $u$ through $\|u\|_{L^{m^*_0}(\Omega)}.$ This way, 
$\m$ is a measure for which Lemma~\ref{lem2.2} could be applied.

For any positive  $k\geq1$ we define now the truncated function
$$
u_k(x)=\max\{u(x)-k,0\}
$$
and the corresponding upper level set
$$
A_k=\{x\in\Omega\colon\ u(x)>k\}.
$$

Noting that boundedness from above of the weak solution to \eqref{1.1} is equivalent to $\m(A_k)=0$ for large enough values of $k,$
we aim now to prove that. For, having
$u_k \in W^{1,m}(\Omega)$ and $u_k\equiv 0$ on $\Omega\setminus A_k,$  the H\"older inequality implies
\begin{align*}
\int_\Omega u_k(x)\, d\m =&\ \int_{A_k} u_k(x)\, d\m\\
 \leq&\ \left(\int_{A_k} d\m\right)^{1-\frac{n-m}{(n-m+\varepsilon_0)m}}\left(\int_{A_k} |u_k(x)|^\frac{(n-m+\varepsilon_0)m}{n-m}\, d\m\right)^\frac{n-m}{(n-m+\varepsilon_0)m},
\end{align*}
whence
\begin{align}\label{3.4}
\int_\Omega u_k(x)\, d\m \leq C (\m(A_k))^{1-\frac{n-m}{(n-m+\varepsilon_0)m}}&\left(\left(\int_{A_k}|Du_k(x)|^m\, dx\right)^{1/m}\right.\\
\nonumber &\quad+\left.\left(\int_{A_k}|u_k(x)|^\sigma\, dx\right)^{1/\sigma}\right)
\end{align}
in view of the Adams trace inequality (Lemma~\ref{lem2.2}) with arbitrary $\sigma\in(0,m].$

In order to estimate the $m$-energy of $u_k$ in \eqref{3.4}, we note that 
  $Du=Du_k$ a.e. in $A_k$ and using $u_k$ as a test function in \eqref{1.2} yields
\begin{align}\label{3.5}
\int_{A_k}& \ba\big(x,u(x),Du(x)\big)\cdot Du_k(x)\,dx+\int_{A_k} b\big(x,u(x),Du(x)\big)u_k(x) \, dx\\
\nonumber
     &=\int_{\partial\Omega\cap \overline{A}_k} \psi\big(x,u(x)\big)u_k(x)\, d\Gamma_x. 
\end{align}

It follows from the coercivity condition \eqref{1.3} that
\begin{align}\label{3.6}
\int_{A_k} & \ba\big(x,u(x),Du(x)\big)\cdot Du_k(x)\,dx\\
\nonumber
&\qquad \geq \gamma \int_{A_k}|Du_k(x)|^m\,dx-\Lambda\int_{A_k}|u(x)|^{m^*}\,dx
-\Lambda\int_{A_k}\varphi_1(x)^{\frac{m}{m-1}}\,dx,
\end{align}
while 
$$
0<\frac{u(x)-k}{u(x)}<1\quad \textrm{a.e.\ in}\ A_k, 
$$
and thus
\begin{align*}
\big|b(x,u(x),&Du(x))u_k(x)\big|= 
\big|b(x,u(x),Du_k(x))u(x)\big|\frac{u(x)-k}{|u(x)|}\\
&\leq \Lambda\left( \varphi_2(x)|u(x)|+|u(x)|^{m^*}+|Du_k(x)|^{\frac{m(m^*-1)}{m^*}}|u(x)|\right)\\
&\leq \Lambda\left( \varepsilon |Du_k(x)|^m +C(\varepsilon)|u(x)|^{m^*}+
\varphi_2(x)|u(x)|\right)
\end{align*}
for any $\varepsilon>0$ 
as consequence of the controlled growth assumption \eqref{1.6} and the Young inequality.

Therefore,
\begin{align}\label{3.7}
\int_{A_k}  b\big(x,u(x),&Du(x)\big)u_k(x) \, dx \leq \Lambda\varepsilon \int_{A_k} |Du_k(x)|^m\,dx\\
\nonumber
&\qquad  +\Lambda C(\varepsilon)\int_{A_k}|u(x)|^{m^*}\,dx+\Lambda \int_{A_k} \varphi_2(x)|u(x)|\,dx.
\end{align}

Substituting \eqref{3.6} and \eqref{3.7} into \eqref{3.5} and choosing $\varepsilon>0$ small enough bring us to  the  energy inequality
\begin{align}\label{3.9}
\int_{A_k}|Du_k(x)|^m\, dx \leq C &\Bigg(
\underbrace{\int_{A_k}\varphi_1(x)^\frac{m}{m-1} \; dx}_{I_1:=} + \underbrace{\int_{A_k} \varphi_2(x)|u(x)| \, dx}_{I_2:=}
\\
\nonumber
& \quad +\! \underbrace{\int_{A_k}\!\! |u(x)|^{\frac{nm}{n-m}}\, dx}_{I_3:=} +\!
\underbrace{\int_{\partial\Omega\cap \overline{A}_k}\!\! \psi\big(x,u(x)\big)u_k(x)\, d\Gamma_x}_{I_4:=} \Bigg).
\end{align}

\medskip

\textit{Step 1: Estimates of the volume integrals
 $I_1,I_2$ and $I_3.$} These already have been  estimated in \cite{AFPS} but, for reader's convenience and for the sake of completeness, we will present here briefly the corresponding arguments. Thus, remembering the definition of the measure $\m$ and $k\geq1,$
we have immediately
\begin{equation}\label{3.10}
I_1 \leq \m(A_k)\leq k^m \m(A_k)
\end{equation}
whereas
\begin{align*}
I_2 =&\ \int_{A_k}\varphi_2(x)|u(x)-k+k|\, dx \leq \int_{A_k}\varphi_2(x)|u_k(x)|\, dx + k \int_{A_k}\varphi_2(x)\, dx\\
\leq&\ \int_{A_k} |u_k(x)|\,d\m +k\m(A_k). 
\end{align*}
As for the integral on the right-hand side, Young's inequality applied to 
 \eqref{3.4} gives
\begin{align*}
\int_{A_k} |u_k(x)|\,d\m \leq &\ \varepsilon \int_{A_k}\!\!|Du_k(x)|^m\, dx \\
&\quad +C(\varepsilon)\left(\int_{A_k}\!\!|u_k(x)|^\sigma\, dx\right)^{m/\sigma}+C \m(A_k)\big(\m(A_k)\big)^{\frac{\varepsilon_0}{(m-1)(n-m+\varepsilon_0)}}
\end{align*}
with arbitrary $\varepsilon>0$ and $\sigma\in(0,m].$ On the other hand,
\begin{equation}\label{3.11}
\m(A_k)\leq \m(\Omega)
\leq  |\Omega| +C \left(
\|\varphi_1\|_{L^{p_1,\lambda_1}(\Omega)}^{\frac{m}{m-1}}+
\|\varphi_2\|_{L^{p_2,\lambda_2}(\Omega)}+
\|u\|_{L^{m^*_0}(\Omega)}^{\frac{m^2}{n-m}}\right)
\end{equation}
and therefore, $k\geq 1$  yields
\begin{equation}\label{3.12}
I_2 \leq \varepsilon \int_{A_k}|Du_k(x)|^m\, dx + C k^m \m(A_k) + C(\varepsilon)\left(\int_{A_k}|u_k(x)|^\sigma\, dx\right)^{m/\sigma}
\end{equation}
with arbitrary $\varepsilon>0$ and $\sigma\in(0,m].$

For what concerns the term 
 $I_3$ in \eqref{3.9}, we have
\begin{align*}
I_3 =&\ \int_{A_k} |u(x)|^{\frac{nm}{n-m}}\, dx =
\int_{A_k} |u(x)-k+k|^m|u(x)|^{\frac{m^2}{n-m}}\, dx \\
\leq&\ 2^{m-1}\left(\int_{A_k} |u_k(x)|^m|u(x)|^{\frac{m^2}{n-m}}\, dx + k^m \int_{A_k} |u(x)|^{\frac{m^2}{n-m}}\, dx \right) \\
\leq&\ 2^{m-1} \int_{A_k} |u_k(x)|^m|u(x)|^{\frac{m^2}{n-m}}\, dx + 2^{m-1}k^m  \m(A_k).
\end{align*}
Noting that $|u|^{\frac{m^2}{n-m}} \in L^{1,\theta}(\Omega)$ with $\theta = n - m + \frac{mm^*_0(n-m)-nm^2}{m^*_0(n-m)}>n-m$
by means of \eqref{3.2}, we will estimate
the integral on the right-hand side through
the Adams trace inequality. For, choose a number $r'<m,$ sufficiently close to $m,$ with the property
\begin{equation}\label{r'}
n-m<\frac{m}{r'}(n-r')=n-r'+\frac{(n-r')(m-r')}{r'}<\theta
\end{equation}
and apply Lemma~\ref{lem2.1} to get
$|u|^{\frac{m^2}{n-m}} \in L^{1,n-r'+\frac{(n-r')(m-r')}{r'}}(\Omega).$ Thus the second claim of  Lem\-ma~\ref{lem2.2}, applied with $c(x)=|u(x)|^{\frac{m^2}{n-m}},$
and the H\"older inequality infer 
\begin{align*}
\int_{A_k} |u_k(x)|^m |u(x)|^{\frac{m^2}{n-m}}\, dx \leq&\ C\left(\left(\int_{A_k}\!\! |Du_k(x)|^{r'}\, dx\right)^{m/r'}\!\!+\left(\int_{A_k}\!\!|u_k(x)|^\sigma\, dx\right)^{m/\sigma} \right)\\
 \leq&\ C|A_k|^{\frac{m}{r'}-1}\!\!\int_{A_k}\!\!|D{u_k}(x)|^m \, dx + C \left(\int_{A_k}\!\!|u_k(x)|^\sigma\, dx\right)^{m/\sigma}
\end{align*}
with $C$ depending also on $\left\||u|^{\frac{m^2}{n-m}}\right\|_{L^{1,\theta}(\Omega)}$ which is anyway bounded in terms of $\|u\|_{L^{m^*_0}(\Omega)}$ as it follows from Lemma~\ref{lem2.4} and \eqref{3.2}. Therefore,
$$
I_3 \leq C\left(|A_k|^{\frac{m}{r'}-1}\int_{A_k}|Du_k(x)|^m \, dx + k^m  \m(A_k)+
\left(\int_{A_k}|u_k(x)|^\sigma\, dx\right)^{m/\sigma} \right)
$$
and, remembering \eqref{3.10} and \eqref{3.12}, we get
\begin{align}\label{3.13}
I_1+I_2+I_3 \leq&\ C\left(\left(|A_k|^{\frac{m}{r'}-1}+\varepsilon\right)\int_{A_k}|Du_k(x)|^m \, dx\right.\\
\nonumber
&\qquad + \left. k^m \m(A_k)+\left(\int_{A_k}|u_k(x)|^\sigma\, dx\right)^{m/\sigma} \right) 
\end{align}
with arbitrary $\varepsilon>0$ and for all $k\geq 1$ and all $\sigma\in(0,m].$

\medskip

\textit{Step 2: Estimate of the surface integral $I_4$ in \eqref{3.9}.}
We will estimate it by means of the multiplicative trace inequality, Lemma~\ref{trace}. Recalling that
 $0<\frac{u(x)-k}{u(x)}<1$ a.e. in $A_k$ we use \eqref{1.8} to get that for a.a. $x\in \partial\Omega\cap \overline{A}_k$ it holds
$$
\psi\big(x,u(x)\big)u_k(x)\leq
\big|\psi\big(x,u(x)\big)\big|u(x)\dfrac{u(x)-k}{u(x)}\leq \psi_1(x)|u(x)|+\psi_2|u(x)|^{1+\beta}
$$
and therefore
\begin{equation}\label{3.8}
I_4\leq
\underbrace{\int_{\partial\Omega\cap \overline{A}_k}
\psi_1(x)|u(x)|\,d\Gamma_x}_{I_{4a}:=}+
\underbrace{\int_{\partial\Omega\cap \overline{A}_k}
\psi_2(x)|u(x)|^{1+\beta}\,d\Gamma_x}_{I_{4b}:=}. 
\end{equation}
We have $u(x)\leq |u(x)|^m$ a.e. in $A_k$
since $u(x)>k\geq 1$ there.  Consider now  the measure $d\m_1:=\psi_1(x)\,d\Gamma_x,$ supported on $\partial\Omega$ and defined by
$
\m_1(E)=\int_{E\cap \partial\Omega} \psi_1(x)\,d\Gamma_x$ for any measurable set  $E.$ Thus
$$
I_{4a}\leq 
\int_{\partial\Omega\cap \overline{A}_k}
\psi_1(x)|u(x)|^m\,d\Gamma_x=
\int_{\partial\Omega\cap \overline{A}_k}
|u(x)|^m\,d\m_1
$$
and to estimate the $L^m(\m_1,\partial\Omega\cap \overline{A}_k)$-norm of $u,$ we will  use the multiplicative Sobolev trace embedding as given in Lemma~\ref{trace}. To check its applicability, consider an arbitrary ball 
$B(x,\rho)$ with center at $x\in\R^n$ and of radius $\rho,$ and note that $B(x,\rho)\cap\partial\Omega \subseteq B(x',\rho)\cap\partial\Omega$ with some $x'\in\partial\Omega.$ Direct computations, based on what has been said in Section~\ref{sec2} about Morrey spaces on $\partial\Omega,$ yield
\begin{align*}
\m_1(B(x,\rho))=&\ \int_{B(x,\rho)\cap\partial\Omega}\psi_1(y)\,d\Gamma_y\leq 
\int_{B(x',\rho)\cap\partial\Omega}\psi_1(y)\,d\Gamma_y\\
\leq&\ C
\rho^{\mu_1/q_1}\|\psi_1\|_{L^{q_1,\mu_1}(\partial\Omega)}\left(
\int_{B(x',\rho)\cap\partial\Omega}d\Gamma_y\right)^{1-1/q_1}\\
\leq&\ C\|\psi_1\|_{L^{q_1,\mu_1}(\partial\Omega)}\rho^{\mu_1/q_1+(n-1)(1-1/q_1)},
\end{align*}
where the constant $C$ depends on the Lipschitz norms of the diffeormorphisms locally flattening $\partial\Omega,$ but is independent of $x$ and $\rho.$

We have $\mu_1/q_1+(n-1)(1-1/q_1)>n-m$ since 
$q_1>\frac{n-1-\mu_1}{m-1}$ by means of \eqref{Q}, and setting $\chi_{A_k}$ for the characteristic function of the set $A_k,$ Lemma~\ref{trace} gives
\begin{align*}
\|u\|_{L^m(\m_1,\partial\Omega\cap \overline{A}_k)}=&\ \|\chi_{A_k}u\|_{L^m(\m_1,\partial\Omega)}\\
\leq&\ C\|\chi_{A_k}u\|^{\tau}_{W^{1,m}(\Omega)}
\|\chi_{A_k}u\|^{1-\tau}_{L^{m}(\Omega)}\\
\leq&\ C\|u\|^{\tau}_{W^{1,m}(A_k)}
\|u\|^{1-\tau}_{L^{m}(A_k)}
\end{align*}
with $\tau=\frac{n}{m}-\frac{\mu_1/q_1+(n-1)(1-1/q_1)}{m}<1.$ We apply now the Young inequality together with $(a+b)^{1/m}\leq a^{1/m}+b^{1/m}$ in order to get
\begin{align}\label{2}
&\|u\|_{L^m(\m_1,\partial\Omega\cap \overline{A}_k)}\leq \frac{\varepsilon}{2} \|u\|_{W^{1,m}(A_k)}+\frac{C}{\varepsilon^{\tau/(1-\tau)}} \|u\|_{L^{m}(A_k)}\\
\nonumber
&\qquad \leq \frac{\varepsilon}{2}\left(\int_{A_k}|Du(x)|^mdx\right)^{1/m}+ \left(\frac{\varepsilon}{2}+
\frac{C}{\varepsilon^{\tau/(1-\tau)}}\right)
\left(\int_{A_k}|u(x)|^mdx\right)^{1/m}
\end{align}
for any $\varepsilon>0,$ and the second term in the right-hand side will be interpolated via the Gagliardo--Nirenberg  inequality. Thus, Lemma~\ref{G-N}, with $r=m$ and 
$\langle\langle u\rangle\rangle_\sigma:=\left(\int_{A_k}|u(x)|^\sigma dx\right)^{1/\sigma},$ reads as
$$
\langle\langle u\rangle\rangle_m\leq C'
\Big(\langle\langle Du\rangle\rangle_m+
\langle\langle u\rangle\rangle_\sigma\Big)^\delta \langle\langle u\rangle\rangle_\sigma^{1-\delta}
$$
for each $\sigma\in(0,m]$ with $\delta=\frac{n(m-\sigma)}{mn-\sigma(n-m)}\in (0,1),$
and therefore 
\begin{align*}
&\left(\frac{\varepsilon}{2}+
\frac{C}{\varepsilon^{\tau/(1-\tau)}}\right)
\langle\langle u\rangle\rangle_m\\
&\qquad\quad \leq 
\left(\frac{\varepsilon}{2}\right)^\delta\Big(\langle\langle Du\rangle\rangle_m+
\langle\langle u\rangle\rangle_\sigma\Big)^\delta 
\left(\frac{2}{\varepsilon}\right)^\delta C' \left(\frac{\varepsilon}{2}+
\frac{C}{\varepsilon^{\tau/(1-\tau)}}\right)
\langle\langle u\rangle\rangle_\sigma^{1-\delta}\\
&\qquad\quad \leq \frac{\varepsilon}{2}
\langle\langle Du\rangle\rangle_m\!+\!
\left[\frac{\varepsilon}{2}+C'^{1/(1-\delta)}
\left(\frac{\varepsilon}{2}+
\frac{C}{\varepsilon^{\tau/(1-\tau)}}\right)^{1/(1-\delta)}\left(\frac{2}{\varepsilon}\right)^{\delta/(1-\delta))} \right]
\langle\langle u\rangle\rangle_\sigma
\end{align*}
by the Young inequality. This way,  \eqref{2} yields 
$$
I_{4a}\leq\|u\|^m_{L^m(\m_1,\partial\Omega\cap \overline{A}_k)}\leq
 \varepsilon \int_{A_k} |Du|^mdx+C(\varepsilon)
\left(\int_{A_k}|u(x)|^\sigma dx\right)^{m/\sigma}
$$
for any level $k\geq1,$ any $\varepsilon>0$ and any $\sigma\in(0,m].$

\medskip

A similar procedure applies also to the integral $I_{4b}$ in \eqref{3.8}. Namely,
$$
I_{4b}=\int_{\partial\Omega\cap \overline{A}_k}
\psi_2(x)|u(x)|^{1+\beta}\,d\Gamma_x
=
\int_{\partial\Omega\cap \overline{A}_k}
|u(x)|^{1+\beta}\,d\m_2,
$$
where the measure $d\m_2:=\psi_2(x)\,d\Gamma_x$ is supported on $\partial\Omega$ and given by $\m_2(E)=\int_{E\cap \partial\Omega} \psi_2(x)\,d\Gamma_x.$ We have
\begin{align*}
\m_2(B(x,\rho))=&\ \int_{B(x,\rho)\cap\partial\Omega}\psi_2(y)\,d\Gamma_y\\
\leq&\ C\|\psi_2\|_{L^{q_2,\mu_2}(\partial\Omega)}\rho^{\mu_2/q_2+(n-1)(1-1/q_2)}
\end{align*}
as above, with a constant depending on the Lipschitz regularity of $\partial\Omega,$ but not on  $x$ and $\rho.$ Moreover, \eqref{Q} ensures that 
$\mu_2/q_2+(n-1)(1-1/q_2)>n-m$ and thus applicability of Lemma~\ref{trace} which reduces the estimate of $I_{4b}$ to that of 
$\|u\|^{1+\beta}_{L^{1+\beta}(\m_2,\partial\Omega\cap \overline{A}_k)}.$ For,
we will distinguish between two cases depending on the values of $\beta$ and $m:$

\medskip

\textsc{Case 1:} $1+\beta\leq m.$ We have $u(x)\geq k\geq 1$ on $\partial\Omega\cap \overline{A}_k$ and therefore $|u(x)|^{1+\beta}\leq |u(x)|^m,$ whence
$$
\int_{\partial\Omega\cap \overline{A}_k}
|u(x)|^{1+\beta}\,d\m_2 \leq \int_{\partial\Omega\cap \overline{A}_k}
|u(x)|^{m}\,d\m_2.
$$
The last integral estimates in absolutely the same manner as $I_{4a}$ above, leading thus to
\begin{equation}\label{Est-I_4b_1}
I_{4b}=\|u\|^{1+\beta}_{L^{1+\beta}(\m_2,\partial\Omega\cap \overline{A}_k)}\leq
 \varepsilon \int_{A_k} |Du|^mdx+C(\varepsilon)
\left(\int_{A_k}|u(x)|^\sigma dx\right)^{m/\sigma}
\end{equation}
for any level $k\geq1,$ any $\varepsilon>0$ and any $\sigma\in(0,m].$

\medskip

\textsc{Case 2:} $1+\beta> m.$ We have
$q_2>\frac{m(n-1)}{n(m-1)-\beta(n-m)}>\frac{m(n-1-\mu_2)}{n(m-1)-\beta(n-m)}$
in view of  \eqref{Q} and Lemma~\ref{trace} gives
\begin{equation}\label{4}
\|u\|_{L^{1+\beta}(\m_2,\partial\Omega\cap \overline{A}_k)}
\leq C\|u\|^{\tau}_{W^{1,m}(A_k)}
\|u\|^{1-\tau}_{L^{m}(A_k)},
\end{equation}
with $\tau=\frac{n}{m}-\frac{\mu_2/q_2+(n-1)(1-1/q_2)}{1+\beta}<1.$ We pick now a real number $\theta$ such that
$$
\theta>\max\left\{\frac{1}{m},\frac{1}{\tau(1+\beta)}\right\},
$$
which additionally satisfies $\theta<\frac{1}{\tau(1+\beta)-(1+\beta-m)}$ whenever $\tau(1+\beta)-(1+\beta-m)>0.$ Actually, that choice of $\theta$ is possible since $1+\beta>m$ and $\tau<1,$ and we have 
$\tau-\frac{1}{\theta(1+\beta)}>0$ and 
$1-\tau-\frac{\theta m-1}{\theta(1+\beta)}>0.$ This allows to decompose
the right-hand side of \eqref{4} as
$$
\|u\|^{\tau}_{W^{1,m}(A_k)}
\|u\|^{1-\tau}_{L^{m}(A_k)}=
\|u\|_{W^{1,m}(A_k)}^{\frac{1}{\theta(1+\beta)}}\|u\|_{L^{m}(A_k)}^{\frac{\theta m-1}{\theta(1+\beta)}}
\|u\|_{W^{1,m}(A_k)}^{\tau-\frac{1}{\theta(1+\beta)}}\|u\|_{L^{m}(A_k)}^{1-\tau-\frac{\theta m-1}{\theta(1+\beta)}},
$$
and setting
$$
N(u):= \|u\|_{W^{1,m}(\Omega)}^{1-\frac{m}{1+\beta}}=
\|u\|_{W^{1,m}(\Omega)}^{\tau-\frac{1}{\theta(1+\beta)}}\|u\|_{W^{1,m}(\Omega)}^{1-\tau-\frac{\theta m-1}{\theta(1+\beta)}},
$$
we get
\begin{align*}
\|u\|_{W^{1,m}(A_k)}^{\tau-\frac{1}{\theta(1+\beta)}}\|u\|_{L^{m}(A_k)}^{1-\tau-\frac{\theta m-1}{\theta(1+\beta)}}\leq &\
\|u\|_{W^{1,m}(A_k)}^{\tau-\frac{1}{\theta(1+\beta)}}\|u\|_{W^{1,m}(A_k)}^{1-\tau-\frac{\theta m-1}{\theta(1+\beta)}}\\
\leq&\ 
\|u\|_{W^{1,m}(\Omega)}^{\tau-\frac{1}{\theta(1+\beta)}}\|u\|_{W^{1,m}(\Omega)}^{1-\tau-\frac{\theta m-1}{\theta(1+\beta)}}= N(u).
\end{align*}
This way \eqref{4} reads
$$
\|u\|_{L^{1+\beta}(\m_2,\partial\Omega\cap \overline{A}_k)}
\leq C N(u)\|u\|_{W^{1,m}(A_k)}^{\frac{1}{\theta(1+\beta)}}\|u\|_{L^{m}(A_k)}^{\frac{\theta m-1}{\theta(1+\beta)}},
$$
whence
\begin{align*}
\int_{\partial\Omega\cap \overline{A}_k}
\psi_2(x)|u(x)|^{1+\beta}d\Gamma_x \leq&\ C
\|u\|_{W^{1,m}(A_k)}^{\frac{1}{\theta}}\|u\|_{L^{m}(A_k)}^{\frac{\theta m-1}{\theta}} N(u)^{1+\beta}\\
\leq&\ \frac{\varepsilon}{2}
\|u\|_{W^{1,m}(A_k)}^{m}+C(\varepsilon,N(u))\|u\|_{L^{m}(A_k)}^{m}\\
\leq&\ \frac{\varepsilon}{2}\int_{A_k}\!\!|Du(x)|^mdx+\left(\frac{\varepsilon}{2}+C(\varepsilon,N(u))\right) \int_{A_k}\!\!|u(x)|^mdx
\end{align*}
by means of the Young inequality and  $\theta m>1.$ This way, interpolating the term $\int_{A_k}|u(x)|^mdx$ by the Gagliardo--Nirenberg inequality (Lemma~\ref{G-N}), we arrive once again to \eqref{Est-I_4b_1} with
the constant $C(\varepsilon)$ depending also on $\|u\|_{W^{1,m}(\Omega)}$ which is not restrictive as it follows from Lemma~\ref{lem2.4}.

Therefore, remembering \eqref{3.8}, the surface integral $I_4$ in \eqref{3.9} satisfies the bound
$$
I_{4}\leq
 \varepsilon \int_{A_k} |Du|^mdx+C(\varepsilon)
\left(\int_{A_k}|u(x)|^\sigma dx\right)^{m/\sigma}
$$
for any level $k\geq1,$ any $\varepsilon>0$ and any $\sigma\in(0,m].$

Indeed, we have $|Du(x)|=|Du_k(x)|$ for a.a. $x\in A_k,$ while
\begin{align*}
\left(\int_{A_k}|u(x)|^\sigma dx\right)^{m/\sigma}=&\ 
\left(\int_{A_k}|u(x)-k+k|^\sigma dx\right)^{m/\sigma}\\
\leq&\  C\left(
\int_{A_k}|u_k(x)|^\sigma dx + k^\sigma\m(A_k)\right)^{m/\sigma}\\
\leq&\ C\left(
\int_{A_k}|u_k(x)|^\sigma dx\right)^{m/\sigma} + Ck^m \m(A_k)^{m/\sigma}\\
\leq&\ C\left(
\int_{A_k}|u_k(x)|^\sigma dx\right)^{m/\sigma} + Ck^m \m(A_k)
\end{align*}
since $\m(A_k)^{m/\sigma}=\m(A_k)\m(A_k)^{m/\sigma-1}\leq \m(A_k)\m(\Omega)^{m/\sigma-1}\leq C\m(A_k).$

Summarizing, we have
\begin{equation}\label{Est-I_4}
I_{4}\leq
 \varepsilon \int_{A_k} |Du|^mdx+C(\varepsilon)
\left(\int_{A_k}|u(x)|^\sigma dx\right)^{m/\sigma}+C(\varepsilon)k^m \m(A_k)
\end{equation}
for any level $k\geq1,$ any $\varepsilon>0$ and any $\sigma\in(0,m].$

\medskip

\textit{Step 3: Decay estimate for the total mass of the weak solution over the upper level set $A_k.$} Putting together \eqref{3.13} and \eqref{Est-I_4} rewrites \eqref{3.9} into
\begin{align}\label{3.13'}
\int_{A_k} |Du_k(x)|^m \, dx
\leq&\ C\left(\left(|A_k|^{\frac{m}{r'}-1}+\varepsilon\right)\int_{A_k}|Du_k(x)|^m \, dx\right.\\
\nonumber
&\qquad + \left. k^m \m(A_k)+\left(\int_{A_k}|u_k(x)|^\sigma\, dx\right)^{m/\sigma} \right) 
\end{align}
with arbitrary $\varepsilon>0,$  for all $k\geq 1$ and all $\sigma\in(0,m],$ where $r'<m$ is choosen according to \eqref{r'}. To manage the quantity $|A_k|^{\frac{m}{r'}-1}$ above, we note that 
\begin{align*}
k^\frac{nm}{n-m}|A_k| \leq&\ \int_{A_k} |u(x)|^\frac{nm}{n-m}\, dx\\
\leq&\ \int_\Omega |u(x)|^\frac{nm}{n-m}\, dx \leq C\Big(\|Du\|_{L^m(\Omega)}+\|u\|_{L^m(\Omega)}\Big)^\frac{nm}{n-m}
\end{align*}
as consequence of the Sobolev inequality (Lemma~\ref{lem2.3}), and this
 means that the term $C\left(|A_k|^{\frac{m}{r'}-1}+\varepsilon\right)$ on the right-hand side of \eqref{3.13'} can be made less than $1/2$ whenever $k \geq k_0$ for  a large enough $k_0,$
depending on known quantities, on $\|u\|_{L^m(\Omega)}$ and on $\|Du\|_{L^m(\Omega)},$ and if $\varepsilon>0$ is small enough.

Thus \eqref{3.13'} becomes
$$
\int_{A_k} |Du_k(x)|^m \, dx \leq C
\left(k^m \m(A_k)+\left(\int_{A_k}|u_k(x)|^\sigma\, dx\right)^{m/\sigma} \right)
\qquad\forall\ k \geq k_0
$$
for all $k\geq k_0$
and \eqref{3.4} takes on the form
\begin{align}\label{3.14}
&\int_\Omega u_k(x)\, d\m \\
\nonumber
&\qquad \leq C \big(\m(A_k)\big)^{1-\frac{n-m}{m(n-m+\varepsilon_0)}}\left(k \big(\m(A_k)\big)^{\frac{1}{m}}
+\left(\int_{A_k}|u_k(x)|^\sigma\, dx\right)^{1/\sigma}\right)
\end{align}
 for all $k\geq k_0$ and all $\sigma\in(0,m].$

At this point, we will show that an appropriate choice of $\sigma\in(0,m]$ yields
\begin{equation}\label{3.15}
\left(\int_{A_k}|u_k(x)|^\sigma\, dx\right)^{1/\sigma} \leq C k \big(\m(A_k)\big)^{\frac{1}{m}}
\end{equation}
with a constant $C$ depending on known quantities. For, two different cases will be distinguished:

\medskip

\textsc{Case 1:} $m<n/2.$ With $\sigma=\frac{m^2}{n-m},$ one has
\begin{align*}
\int_{A_k}|u_k(x)|^{\frac{m^2}{n-m}}\,dx =&\ 
	\int_{A_k}|u(x)-k|^{\frac{m^2}{n-m}}\,dx \leq C\left(\int_{A_k}|u(x)|^{\frac{m^2}{n-m}}\,dx+k^{\frac{m^2}{n-m}}|A_k|\right)\\
	\leq&\ C\left(\m(A_k)+k^{\frac{m^2}{n-m}}|A_k|\right)\leq C k^{\frac{m^2}{n-m}}\m(A_k)
\end{align*}
since $|A_k|\leq\m(A_k)$ and $k\geq 1,$ whence
\begin{align*}
\left(\int_{A_k}|u_k(x)|^{\frac{m^2}{n-m}}\,dx\right)^{\frac{n-m}{m^2}} \leq&\ C k\big(\m(A_k)\big)^{\frac{n-m}{m^2}}=
C k\big(\m(A_k)\big)^{\frac{1}{m}}
(\m(A_k)\big)^{\frac{n-m}{m^2}-\frac{1}{m}}\\
\leq&\ C k\big(\m(A_k)\big)^{\frac{1}{m}}
(\m(\Omega)\big)^{\frac{n-m}{m^2}-\frac{1}{m}}
\end{align*}
because of $\frac{n-m}{m^2}-\frac{1}{m}>0.$ This \eqref{3.11} implies \eqref{3.15}.

\medskip

\textsc{Case 2:} $m\geq n/2.$ We choose  $\sigma=1$ now and observe that $\frac{m^2}{n-m}\geq 1.$ Then
$$
\int_{A_k}|u_k(x)|\,dx =
	\int_{A_k}|u(x)-k|\,dx \leq \int_{A_k}|u(x)|\,dx+k|A_k|
	\leq \int_{A_k}|u(x)|\,dx + k\m(A_k),
$$
while the H\"older inequality and $\frac{m^2}{n-m}\geq 1$ yield
$$
\int_{A_k}|u(x)|\,dx\leq \left(
\int_{A_k}|u_k(x)|^{\frac{m^2}{n-m}}\,dx\right)^{\frac{n-m}{m^2}}|A_k|^{1-\frac{n-m}{m^2}}\leq \m(A_k).
$$
Since $k\geq 1$ anyway, we get
\begin{align*}
\int_{A_k}|u_k(x)|\,dx \leq&\ 2k\m(A_k)=2k
\big(\m(A_k)\big)^{\frac{1}{m}}
\big(\m(A_k)\big)^{1-\frac{1}{m}}\\
\leq&\ 2k
\big(\m(A_k)\big)^{\frac{1}{m}}
\big(\m(\Omega)\big)^{1-\frac{1}{m}}\leq 
C k \big(\m(A_k)\big)^{\frac{1}{m}}
\end{align*}
by \eqref{3.11}, and therefore \eqref{3.15} holds true also in the second case.

\medskip

Once having \eqref{3.15}, \eqref{3.14} becomes
\begin{equation}\label{3.16}
 \int_{A_k} u_k(x)\, d\m \leq Ck\big(\m(A_k)\big)^{1+\frac{\varepsilon_0}{m(n-m+\varepsilon_0)}}\qquad \forall\ k \geq k_0, 
\end{equation}
where $\varepsilon_0>0.$

The Cavalieri principle gives
$$
\int_{A_k} u_k(x)\, d\m= \int_{A_k}(u(x)-k)\, d\m=\int_k^\infty \m(A_t)\, dt,
$$
whence, setting $\tau(t):=\m(A_t),$  \eqref{3.16} rewrites into
$$
\int_k^\infty \tau(t)\, dt \leq Ck\tau(k)^{1+\delta}\qquad \forall k\geq k_0
$$
with $\delta=\frac{\varepsilon_0}{m(n-m+\varepsilon_0)}>0.$
It remains to apply the Hartman--Stampacchia maximum principle (Lemma~\ref{lem2.5}) in order to conclude that
$$
u(x) \leq k_{\max}\qquad \mathrm{a.e.}\ \Omega
$$
where $k_{\max}$ depends on known quantities, on $\|u\|_{L^m(\Omega)}$ and on  $\|Du\|_{L^m(\Omega)}$ in addition.

Applying the above procedure to $-u(x)$ instead of $u(x),$ yields a 
 bound from below for $u(x)$ as well,   and this gives the  claim \eqref{3.1} when $m<n.$

\bigskip

The limit case $m=n$ is to be treated in the same manner as already done, slightly changing the approach  adopted. Precisely,
 the coercivity condition \eqref{1.3}
and the controlled growth assumption  \eqref{1.6}  take now the form
\begin{align}\label{3.17}
\ba(x,z,\xi)\cdot\xi\geq&\ \gamma|\xi|^{n}-\Lambda|z|^{{m^*}}-\Lambda\varphi_1(x)^\frac{n}{n-1},\\
\label{3.18}
|b(x,z,\xi)| \leq&\ \Lambda\left(\varphi_2(x)+|z|^{{m^*}-1}+|\xi|^{\frac{n({m^*}-1)}{{m^*}}}\right),
\end{align}
where ${m^*}>n$ is an \textit{arbitrary} exponent,
$\varphi_1\in L^{p_1,\lambda_1}(\Omega)$ with $p_1>\frac{n}{n-1},$ $\lambda_1\in[0,n)$ and  $\varphi_2\in L^{p_2,\lambda_2}(\Omega)$ with $p_2> 1,$  $\lambda_2\in[0,n).$

Without loss of generality we may select a number $m'<n,$ close enough to $n,$ in a way that ${m^*}=\frac{n^2}{(n-m')(n+1)}.$ With this choice of $m'$ one has
$$
{m^*}<{(m')^*}=\frac{nm'}{n-m'},\quad \frac{n}{n-1}<\frac{m'}{m'-1},\quad \frac{n({m^*}-1)}{{m^*}}=\frac{m'({(m')^*}-1)}{(m')^*}
$$
and \eqref{3.17} becomes
\begin{align}\label{3.19}
\ba(x,z,\xi)\cdot\xi\geq&\ \gamma|\xi|^{n}-\Lambda|z|^{{m^*}}-\Lambda\varphi_1(x)^\frac{n}{n-1}\\
\nonumber
   \geq&\ \gamma|\xi|^{m'}-\Lambda|z|^{(m')^*}-\Lambda\varphi_1(x)^\frac{m'}{m'-1}
\end{align}
when $|z|\geq 1$ and $|\xi|\geq1$ and where, without loss of generality, we have supposed $\varphi_1(x)\geq1.$ Similarly,
 \eqref{3.18} reads
\begin{equation}\label{3.20}
|b(x,z,\xi)|  \leq \Lambda\left(\varphi_2(x)+|z|^{(m')^*-1}+|\xi|^{\frac{m'((m')^*-1)}{(m')^*}}\right)
\end{equation}
whenever $|z|\geq 1$ and $|\xi|\geq1.$  

Considering the measure
$$
d\m'=\left(\chi_\Omega(x)+\varphi_1(x)^\frac{m'}{m'-1}+\varphi_2(x)+|u(x)|^\frac{m'^2}{n-m'}\right)\,dx,
$$
we may increase eventually the value of $m',$ maintaining  it anyway strictly less than $n,$ in order to have
$p_1>\frac{m'}{m'-1},$ $(m'-1)p_1+\lambda_1>n$ and $m'p_2+\lambda_2>n.$ This leads to
$$
\m'(B_\rho)\leq K \rho^{n-m'+\varepsilon_0}
$$
as before with a suitable $\varepsilon_0>0.$
Defining the functions $u_k(x)$ and the sets $A_k$ as in the case $m<n,$
 it is clear that
$$
\int_{\{x\in A_k\colon |Du_k(x)|<1\}} |Du_k(x)|^{m'}\,dx\leq |A_k|\leq k^{m'}\m'(A_k),
$$
while 
$$
\int_{\{x\in A_k\colon |Du_k(x)|\geq1\}} |Du_k(x)|^{m'}\,dx
$$
can be estimated with the help of \eqref{3.19} and \eqref{3.20}, as already done before, leading this way to \eqref{3.9} with $m'$ instead of $m.$ Regarding the 
volume integrals $I_1,$ $I_2$ and $I_3$ there, the procedure applied in the case $m<n$ yields a bound, 
essentially of the same type as \eqref{3.13} with $m'$ in the place of $m.$

In order to estimate the surface integral $I_4,$ we note that $\beta\in[0,+\infty)$ in \eqref{1.8} when $m=n,$ where $\psi_i\in L^{q_i,\mu_i}(\partial\Omega)$ with $q_i>1$ and $\mu_i\in[0,n-1),$ $i=1,2$ 
(cf. \eqref{Q}). 
Applying the approach already used to get
\eqref{Est-I_4} gives a bound of the same type with $m'$ instead of $m$ increasing, if necessary,  the value of $m'<n$ in a way to have $\mu_i/q_i+(n-1)/(1-1/q_i)>n-m'.$

Definitely, also when $m=n$ we arrive at the estimate  \eqref{3.16} with $m'$ instead of $m,$ and the Hartman--Stampacchia maximum principle Lemma~\ref{lem2.5} leads to the desired bound \eqref{3.1}. This completes the proof of Theorem~\ref{thm3.1}.

\section*{Acknowledgements}

\textit{Dian Palagachev} and \textit{Lubomira  Softova} are members of INdAM-GNAMPA. 

The work of 
\textit{Dian Palagachev} 
is partially supported by the Italian Ministry of Education, University and Research
under the Programme ``Department of Excellence'' L. 232/2016 (Grant No. CUP -
D93C23000100001), by the Research Project of National Relevance ``Evolution problems involving interacting scales'' granted by the Italian Ministry of Education, University and Research (MIUR Prin 2022, project code 2022M9BKBC, Grant No. CUP  D53D23005880006), and by the National Recovery and Resilience Plan (NRRP) funded by the European Union -- NextGenerationEU --
Project Title ``Mathematical Modeling of Biodiversity in the Mediterranean sea: from bacteria to predators, from meadows to currents'' -- project code P202254HT8 -- CUP B53D23027760001 --
Grant Assignment Decree No. 1379 adopted on 01/09/2023 by the Italian Ministry of University and Research (MUR).

The  research of \textit{Lubomira  Softova} is partially supported  by the GNAMPA 2025 Project
``Regularity of the solutions of parabolic equations with non-standard degeneracy''
(Grant No. CUP - E5324001950001).

\end{document}